\documentclass{amsart}

\usepackage{amsmath,amssymb,amsfonts,enumerate,amsthm}

\DeclareMathOperator{\hight}{ht}

\newcommand{\id}{\mbox{id}\,}
\newcommand{\fd}{\mbox{fd}\,}
\newcommand{\rfd}{\mbox{Rfd}\,}

\newcommand{\gfd}{\mbox{Gfd}\,}
\newcommand{\gid}{\mbox{Gid}\,}
\newcommand{\CM}{\mbox{CM}\,}
\newcommand{\CMfd}{\mbox{CMfd}\,}
\newcommand{\CIfd}{\mbox{CIfd}\,}
\newcommand{\CMid}{\mbox{CMid}\,}
\newcommand{\CIid}{\mbox{CIid}\,}
\newcommand{\pd}{\mbox{pd}\,}

\newcommand{\cmdim}{\mbox{CM-dim}\,}
\newcommand{\cmd}{\mbox{cmd}\,}

\newcommand{\gd}{\mbox{G-dim}\,}
\newcommand{\cid}{\mbox{CI-dim}\,}

\newcommand{\Hom}{\mbox{Hom}\,}
\newcommand{\Tor}{\mbox{Tor}\,}
\newcommand{\Ext}{\mbox{Ext}\,}

\newcommand{\Ass}{\mbox{Ass}\,}

\newcommand{\Spec}{\mbox{Spec}\,}
\newcommand{\Supp}{\mbox{Supp}\,}
\newcommand{\depth}{\mbox{depth}\,}
\newcommand{\width}{\mbox{width}\,}
\newcommand{\grade}{\mbox{grade}\,}
\renewcommand{\dim}{\mbox{dim}\,}

\renewcommand{\H}{\mbox{H}}

\newcommand{\R}{\mathbb{R}}

\newcommand{\fn}{\frak{n}}
\newcommand{\fm}{\frak{m}}
\newcommand{\fp}{\frak{p}}

\newcommand{\fq}{\frak{q}}
\newtheorem{thm}{Theorem}[section]
\newtheorem{cor}[thm]{Corollary}
\newtheorem{lem}[thm]{Lemma}
\newtheorem{prop}[thm]{Proposition}
\newtheorem{nota}[thm]{Notation}
\newtheorem{defn}[thm]{Definition}

\newtheorem{exam}[thm]{Example}
\newtheorem{rem}[thm]{Remark}
\newtheorem{ques}[thm]{Question}

\begin{document}

\bibliographystyle{amsplain}

\date{}

\author{Parviz Sahandi, Tirdad Sharif, and Siamak Yassemi}

\address{Department of Mathematics, University of
Tehran\\ P.O. Box 13145--448, Tehran, Iran.}

\email{sahandi@ipm.ir}

\address{Institute for
studies in Theoretical Physics and Mathematics (IPM), Tehran Iran.}

\email{sharif@ipm.ir}

\address{Department of Mathematics, University of
Tehran\\ P.O. Box 13145--448, Tehran, Iran and Institute for studies
in Theoretical Physics and Mathematics (IPM), Tehran Iran.}

\email{yassemi@ipm.ir}



\thanks{T. Sharif was supported in part by a grant from IPM (No.
83130311)}
\thanks{S. Yassemi was supported in part by a grant from
IPM (No. 861300000).}
\title{homological flat dimensions}

\begin{abstract}

For finitely generated module $M$ over a local ring $R$, the
conventional notions of complete intersection dimension $\cid_R M$
and Cohen-Macaulay dimension $\cmdim_R M$ do not extend to cover the
case of infinitely generated modules. In this paper we introduce
similar invariants for not necessarily finitely generated modules,
(namely, complete intersection flat and Cohen-Macaulay flat
dimensions) which for finitely generated modules, coincide with the
corresponding classical ones.

\end{abstract}

\maketitle

\section{Introduction}

An important motivation for studying homological dimensions goes
back to 1956 when Auslander, Buchsbaum and Serre proved the
following theorem: A commutative noetherian local ring $R$ is
regular if the residue field $k$ has finite projective dimension
and only if all $R$-modules have finite projective dimension. This
introduced the theme that finiteness of a homological dimension
for all modules singles out rings with special properties.

Auslander and Bridger \cite{AB}, introduced a homological dimension
designed to single out modules with properties similar to those of
modules over Gorenstein rings.  They called it $\mbox{G}$-dimension
and it is a refinement  of the projective dimension and showed that
a local noetherian ring $(R,\fm,k)$ is Gorenstein if the residue
field $k$ has finite $\mbox{G}$-dimension and only if all finitely
generated $R$-modules have finite $\mbox{G}$-dimension. More
recently, other homological dimensions have been introduced  to
characterize complete intersection and Cohen-Macaulay rings (see
\cite{AGP}, \cite{G} and \cite{A1} for an overview.)

This paper is concerned with homological dimensions for not
necessarily \emph{finitely generated} modules over commutative
noetherian local rings $(R,\fm,k)$ with identity. For any $R$-module
$M$, the flat dimension of $M$ over $R$ is denoted by $\fd_R M$.
There is always an inequality $\fd_R M\leqslant\pd_R M$, and
equality holds if $M$ is \emph{finite}, that is finitely generated,
where $\pd_RM$ denotes for projective dimension of $M$. A deep
result, due to Gruson–-Raynaud \cite{RG} and Jensen \cite{J}, says
that flat R-modules have finite projective dimension. Hence the flat
dimension and the projective dimension of a module are finite
simultaneously. Therefore it seems that, the flat dimension is a
good and suitable extension of the projective dimension for
non-finite modules.

In \cite{CFF} Christensen, Foxby, and Frankild introduced the
\emph{large restricted flat dimension} which is denoted by $\rfd$
and it is defined by the formula
$$\rfd_R M = \sup\{i| \Tor_i ^R(L,M) \neq 0\ \makebox{\rm for some
$R$-module $L$ with $\fd_R L<\infty$} \}.$$ They showed that for all
$R$-module $M$, there is an inequality $$\rfd_RM\le\fd_RM$$ with
equality if $\fd_RM<\infty$.

In \cite {E} and \cite{EJT} Enochs and Jenda have introduced the
Gorenstein flat dimension $\gfd_R M$ of any R-module M. An
$R$-module $M$ is said to be Gorenstein flat if and only if there is
an exact sequence
$$\cdots \rightarrow F^{-1}\rightarrow F^0\rightarrow F^1\rightarrow
\cdots$$ of flat $R$-modules such that $M=\ker(F^0\rightarrow
F^1)$ and such that for any injective $R$-module $I$, $I\otimes _R
-$ preserves the exactness of the above complex. The Gorenstein
flat dimension is defined by using Gorenstein flat modules in a
fashion similar to that of flat dimension.  Recall that for a
finite $R$-module $M$ we have $\gfd_RM=\gd_RM$ by \cite{EJT}. Holm
has studied this concept further in \cite {H} and proved that
$\gfd_R M$ is a refinement of $\fd_R M$ and that $\rfd_R M$ is a
refinement of $\gfd_R M$.  In other words, for any R-module M
there is a chain of inequalities $$\rfd_R M\le\gfd_R M\le\fd_R
M,$$ and if one of these quantities is finite then there is
equality everywhere to its left.

The main goal of this paper is to introduce and study notions of
complete intersection flat dimension ($\CIfd$) and Cohen-Macaulay
flat dimension ($\CMfd$) as refinements of flat dimensions for every
module $M$ over a noetherian ring $R$ (see Sections 3 and 4 for
definitions).

A main result of this paper is the comparison of the Gorenstein
flat and the complete intersection flat dimensions as given by the
following theorem (see Theorem \ref{cmi}):

\vspace{.05in} \noindent{\bf Theorem A.} Let $M$ be an $R$-module.
Then there is an inequality
$$\gfd_RM\le\mbox{CI}\fd_RM$$ with equality if $\mbox{CI}\fd_RM$ is
finite.

\noindent Viewing the above theorem, there is the following sequence
of inequalities:
$$\rfd_RM\le\CMfd_RM\le\gfd_RM\le\CIfd_RM\le\fd_RM.$$
If one of these dimensions is finite, then it is equal to those of
its left.

We also introduce and study a variety of refinements of flat
dimension, namely upper Cohen-Macaulay flat dimension ($\CM^*\fd$)
and upper Gorentein flat dimension ($\mbox{G}^*\fd$) for every
module $M$ over a noetherian ring $R$ (see Sections 3 for
definitions). These dimensions fit into the following scheme of
inequalities:
$$\rfd_RM\le\CM^*\fd_RM\le\mbox{G}^*\fd_RM\le\mbox{CI}\fd_RM\le\fd_RM,$$
with equality to the left of any finite number.

The new homological flat dimensions are in many respects, similar to
the classical ones. As a second example of what can be gained from
our homological flat dimensions, we have the following result which
is called Intersection Theorem for homological flat dimensions; (see
Theorem \ref{int}):

\vspace{.05in} \noindent{\bf Theorem B.} Let $M$ be an $R$-module,
with $\H\fd_R M<\infty$ and of finite depth. Suppose that $R$ is an
equicharacteristic zero ring, then:
$$\dim R \le\dim_R M+ \H\fd_R M,$$
for $\H=\mbox{CI}$, $\mbox{G}^*$, and $\CM^*$.

In Section 4, a number of base change results for homological flat
dimensions are obtained (see Propositions \ref{P1}, \ref{P3}).
Special attention is given to finite homomorphisms $R\to
R/(x_1,\cdots,x_n)$ where $x_1,\cdots,x_n$ is an $R$-regular
elements (see Propositions \ref{P4} and \ref{P5}).

The Auslander-Buchsbaum formula asserts that if a finitely
generated $R$-module $M$ has finite projective dimension, then
$\depth_RM+\pd_RM=\depth R$. In \cite{A2} Auslander further
generalized this formula for $M$ as before and $N$ a finitely
generated $R$-module.  In fact he showed that for
$s=\sup\{n|\Tor_n^R(M,N)\neq0\}$ if either $s=0$ or
$\depth_R\Tor_s^R(M,N)\le 1$, then
$$(*) \qquad s=\depth R-\depth_RM-\depth_RN+\depth_R\Tor_n^R(M,N).$$
More generally we say that the \emph{depth formula} holds for $M$
and $N$ if $s$ is finite and $(*)$ holds. In Section 5 the following
result which may be regarded as an analogue of Auslander's theorem
for $\CIfd$ is proven; (see Theorem \ref{D2}):

\vspace{.05in} \noindent{\bf Theorem C.} Let $M$ and $N$ be
$R$-modules such that $\CIfd_R M<\infty$. If $s$ is finite, then
$$s\geq\depth R-\depth_R M-\depth_R N$$ with equality if and only
if $\depth_R \Tor_s ^R (M,N)=0$.

This is an extension (to non-finite case) of \cite[Theorem
(2.2)]{JJ}.

In Section 6 basic properties of homological flat dimensions for
finitely generated modules are established and in Section 7 we
discuss various homological injective dimensions for modules over
noetherian rings namely, the Cohen-Macaulay injective dimension
($\CMid$), upper Cohen-Macaulay injective dimension ($\CM^*\id$),
upper Gorentein injective dimension ($\mbox{G}^*\id$) and complete
intersection injective dimension ($\mbox{CI}\id$).  These dimensions
satisfy the following inequalities
$$\mbox{Ch}\id_RM\le\CMid_RM\le\CM^*\id_RM\le\mbox{G}^*\id_RM\le\mbox{CI}\id_RM\le\id_RM,$$
where $$\mbox{Ch}\id_RM=\sup \{\depth R_\fp - \width _{R_{\fp}}
M_\fp | \fp\in \Spec (R)\}.$$ Recall that
$\width_RM=\inf\{i|\Tor_i^R(M,k)\neq 0\}$.

It is natural to ask when homological flat dimensions satisfy a
formula of Auslander-Buchsbaum type. The answer is given in the
following theorem; (see Theorem \ref{AB}):

\vspace{.05in} \noindent{\bf Theorem D.} Let $R$ be a Cohen-Macaulay
local ring and let $M$ be an $R$-module of finite $\H\fd_RM$ for
$\H=\mbox{CI}$, $\mbox{G}^*$, $\CM^*$, and $\CM$. Then
$\H\fd_RM+\depth_RM=\depth R$ if and only if
$\depth_RM\le\grade(\fp,M)+\dim R/\fp$ for all $\fp\in\Supp(M)$.

\section{Definitions and Notations}

In this section we recall various definitions of homological
dimensions for finite modules.

\begin{defn} A finite $R$-module $M$ has $G$-dimension 0 if the following
conditions are satisfied:
\begin{itemize}
\item[(i)] $M\cong \Hom_R (\Hom_R (M,R),R)$,

\item[(ii)] $\Ext_R ^i (M,R)=0$ for all $i>0$, and

\item[(iii)] $\Ext_R ^i (\Hom_R(M,R),R)=0$ for all $i>0$.
\end{itemize}
\end{defn}
The Gorenstein dimension of $M$ which is defined by Auslander and
Bridger \cite{AB} and denoted by $\gd_RM$, as the least number $n$
for which there exists an exact sequence
$$0\rightarrow G_n\rightarrow G_{n-1}\rightarrow \cdots
\rightarrow G_0 \rightarrow M\rightarrow 0,$$ where $G_i$ has
$G$-dimension 0 for $i=0,\cdots,n$.

A finite $R$-module $M$ is called \textit{perfect} (resp.
\textit{G-perfect}) if $\pd_R M=\grade_R M$ (resp. $\gd_R
M=\grade_R M$). Let $Q$ be a local ring and $J$ an ideal of $Q$.
By abuse of language we say that $J$ is \textit{perfect} (resp.
\textit{G-perfect}) if the $Q$-module $Q/J$ has the corresponding
property.

The ideal $J$ is called Gorenstein if it is perfect and $\beta_g
^Q (Q/J)=1$ for $g=\grade_Q J$, where $\beta_g ^Q (Q/J)$ stands
for $g$-th betti number of $Q/J$. It is called complete
intersection ideal, if $J$ is generated by an $R$-regular
elements.

We say that $R$ has a $\mbox{CI}$-deformation (resp.
$\mbox{G}^*$-deformation, $\CM$-deformation) if there exists a
local ring $Q$ and a complete intersection (resp. Gorenstein,
G-perfect) ideal $J$ in $Q$ such that $R=Q/J$. A
$\mbox{CI}$-quasi-deformation (resp.
$\mbox{G}^*$-quasi-deformation, $\CM$-quasi-deformation) of $R$ is
a diagram of local homomorphisms $R\rightarrow R'\leftarrow Q$,
with $R\rightarrow R'$ a flat extension and $R'\leftarrow Q$ a
$\mbox{CI}$-deformation (resp. $\mbox{G}^*$-deformation,
$\CM$-deformation). We set $M'=M\otimes_R R'$.

\noindent The \textit{complete intersection dimension} of $M$ as
defined by Avramov, Gasharov, and Peeva \cite{AGP} and denoted by
$\cid_RM$ is
$$\cid_R M:=\inf\{\pd_Q M' -\pd_Q R'| \text{ }R\rightarrow R'\leftarrow
Q \text{ is a }\mbox{CI}\text{-quasi-deformation}\}.$$

\noindent The \textit{upper Gorenstein dimension} of $M$ as
defined by Veliche \cite{V} and denoted by $\mbox{G}^*
\text{-}\dim_R M$ is
$$\mbox{G}^* \text{-}\dim_R M:=\inf\{\pd_Q M' -\pd_Q R'| \text{ }R\rightarrow R'\leftarrow
Q \text{ is a }\mbox{G}^*\text{-quasi-deformation}\}.$$

\noindent The \textit{Cohen-Macaulay dimension} of $M$, as defined
by Gerko \cite{G} and denoted by $\cmdim_R M$ is
$$\cmdim_R M:=\inf\{\gd_Q M' -\gd_Q R'| \text{ }R\rightarrow R'\leftarrow
Q \text{ is a }\CM\text{-quasi-deformation}\}.$$

\noindent  There are the following sequence of inequalities:
$$\cmdim_R M \le\gd_R M \le\mbox{G}^* \text{-}\dim_R M\le\mbox{CI}\text{-}\dim_R M\le\pd_R M,$$
with equality to the left of any finite number.

In \cite{EJ1} Enochs and Jenda introduced the Gorenstein injective
dimension $\gid_RM$ of any $R$-module $M$ as follows:

\begin{defn}
An $R$-module $M$ is said to be Gorenstein injective if and only
if there is an exact sequence $$\cdots \rightarrow
E^{-1}\rightarrow E^0\rightarrow E^1\rightarrow \cdots$$ of
injective $R$-modules such that $M=\ker(E^0\rightarrow E^1)$ and
for any injective $R$-module $E$, $\Hom_R (E,-)$ preserves
exactness of the above complex. The Gorenstein injective dimension
is defined by using Gorenstein injective modules in a fashion
similar to that of injective dimension.
\end{defn}
\noindent It is known that $\gid_RM\le\id_RM$ with equality if
$\id_RM$ is finite.

For a notherian ring the following categories were introduced by
Avramov and Foxby \cite{AF2}:

\begin{defn} Let $R$ be a ring with a dualizing complex $D$. Let
$\mathcal{D}_b(R)$ denote the full subcategory of $\mathcal{D}(R)$
(the derived category of $R$-complexes) consisting of complexes
$X$ with $\H_n(X)=0$ for $n\gg0$. The \emph{Auslander class}
$\mathbf{A}(R)$ is defined as the full subcategory of
$\mathcal{D}_b(R)$, consisting of those complexes $X$ for which
$D\otimes^{\mathbf{L}}_RX\in\mathcal{D}_b(R)$ and the canonical
morphism
$$\gamma_X:X\to\mathbf{R}\Hom_R(D,D\otimes_R^{\mathbf{L}}X),$$
ia an isomorphism.  The \emph{Bass class} $\mathbf{B}(R)$ is
defined as the full subcategory of $\mathcal{D}_b(R)$, consisting
of those complexes $X$ for which
$\mathbf{R}\Hom_R(D,X)\in\mathcal{D}_b(R)$ and the canonical
morphism
$$\iota_X:D\otimes^{\mathbf{L}}_R\mathbf{R}\Hom_R(D,X)\to X,$$
is an isomorphism.
\end{defn}

\begin{rem}\label{Re} It is proved in \cite[(4.1) and (4.4)]{CFH} that if $R$ admits of a
dualizing complex then for an $R$-module $M$ we have:
\begin{itemize}
\item[(a)] $M\in \mathbf{A}(R)\text{ if and only if
}\gfd_RM<\infty,$ and

\item[(b)] $M\in \mathbf{B}(R)\text{ if and only if }\gid_RM<\infty.$
\end{itemize}
See also \cite{ET} and \cite{ET1} for an interesting extension of
this result.
\end{rem}

\section{Complete intersection flat
dimension}

In this section we introduce \textit{complete intersection flat
dimension}, \textit{upper Gorenstein flat dimension}, and
\textit{upper Cohen-Macaulay flat dimension} for not necessarily
finite $R$-modules, and verify a number of their properties which
are similar to those for the flat dimension.

We say that $R$ has a $\CM^*$-deformation if there exist a local
ring $Q$ and a perfect ideal $J$ in $Q$ such that $R=Q/J$. A
$\CM^*$-quasi-deformation of $R$ is a diagram of local homomorphisms
$R\rightarrow R'\leftarrow Q$ with $R\rightarrow R'$ a flat
extension and $R'\leftarrow Q$ a $\CM^*$-deformation. We set
$M'=M\otimes_R R'$.

\begin{defn} Let $M\neq0$ be an $R$-module. The \textit{complete
intersection flat dimension}, \textit{upper Gorenstein flat
dimension}, and \textit{upper Cohen-Macaulay flat dimension} of
$M$, are defined as:
$$\mbox{CI}\fd_RM:=\inf\{\fd_Q M' -\fd_Q R'| \text{ }R\rightarrow R'\leftarrow
Q \text{ is a }\mbox{CI}\text{-quasi-deformation}\}$$
$$\mbox{G}^*\fd_RM:=\inf\{\fd_Q M' -\fd_Q R'| \text{ }R\rightarrow R'\leftarrow
Q \text{ is a }\mbox{G}^*\text{-quasi-deformation}\}$$
$$\CM^*\fd_RM:=\inf\{\fd_Q M' -\fd_Q R'| \text{ }R\rightarrow R'\leftarrow Q
\text{ is a }\CM^*\text{-quasi-deformation}\},$$ respectively.  We
complement this by $\H\fd_R 0=-\infty$ for $\H=\mbox{CI}$,
$\mbox{G}^*$, and $\CM^*$.
\end{defn}

Our first result says that the large restricted flat dimension is a
refinement of the above $\H$-flat dimensions.

\begin{prop} \label{R1} Let $R\rightarrow S\leftarrow Q$ be a $\CM$-quasi-deformation,
and let $M$ be an $R$-module. Then
$$\rfd_Q (M\otimes_R S)-\rfd_Q S=\rfd_R M.$$
\end{prop}
\begin{proof} First we prove the equality $$\rfd_SN+\gd_Q S=
\rfd_QN,$$ for an $S$-module $N$.  To this end, choose by \cite
[(2.4)(b)]{CFF} a prime ideal $\fp$ of $S$ such that the first
equality below holds. Let $\fq$ be the inverse image of $\fp$ in
$Q$. Therefore there is an isomorphism $N_{\fp}\cong N_{\fq}$ of
$Q_{\fq}$-modules and a $\CM$-deformation $Q_{\fq}\rightarrow
S_{\fp}$.  Hence
\begin{align*}
\rfd_SN= & \depth S_{\fp}-\depth_{S_{\fp}}N_{\fp} \\[1ex]
        = & \depth_{Q_{\fq}} S_{\fp}-\depth_{Q_{\fq}}N_{\fp} \\[1ex]
        = & \depth Q_{\fq}-\gd_{Q_{\fq}} S_{\fp}-\depth_{Q_{\fq}}N_{\fp} \\[1ex]
        \le & \rfd_QN-\gd_{Q_{\fq}} S_{\fp} \\[1ex]
        = & \rfd_QN-\gd_Q S.
\end{align*}
The second equality holds since $Q_{\fq}\rightarrow S_{\fp}$ is
surjective; the third equality holds by Auslander-Bridger formula
\cite{AB}; the fourth equality is due to the
$\mbox{G}$-perfectness assumption of $S$ over $Q$; while the
inequality follows from \cite [(2.4)(b)]{CFF}. Now by
\cite[(3.5)]{SY2} we have
$$\rfd_QN\le\rfd_SN+\rfd_Q S\le\rfd_QN-\gd_Q S+\rfd_Q S=\rfd_QN,$$
which is the desired equality.

Now we have
\begin{align*}
\rfd_Q (M\otimes_R S)\le & \rfd_S (M\otimes_R S)+\rfd_Q S \\[1ex]
                     = & \rfd_S (M\otimes_R S)+\gd_Q S \\[1ex]
                     = & \rfd_Q (M\otimes_R S),
\end{align*}
where the inequality was proven in \cite [(3.5)]{SY2}, the first
equality follows from the hypotheses, and the last follows from the
above observation.  Hence $$\rfd_Q (M\otimes_R S)-\rfd_Q S=\rfd_S
(M\otimes_R S)=\rfd_RM$$ where the second equality holds by
\cite[(8.5)]{IS}.
\end{proof}

\begin{thm} \label{CH} Let $M$ be an $R$-module. Then we have $\rfd_RM\le\H\fd_R M$, for $\H=\mbox{CI}$,
$\mbox{G}^*$, and $\CM^*$, with equality if $\H\fd_R M$ is finite.
In this case we have
$$\H\fd_R M=\sup \{\depth  R_\fp - \depth _{R_{\fp}} M_\fp | \fp\in
\Spec (R)\}.$$
\end{thm}
\begin{proof} The inequality follows easily from the definitions of various $\H$-flat dimensions introduced above,
and Proposition
\ref{R1}. The last equality follows from \cite[(2.4)(b)]{CFF}.
\end{proof}

The following corollary is an immediate consequence of Theorem
\ref{CH}:

\begin{cor} \label{IN} There is the following chain of
inequalities:
$$\rfd_RM\le\CM^*\fd_RM\le\mbox{G}^*\fd_RM\le\mbox{CI}\fd_RM\le\fd_RM,$$
with equality to the left of any finite number.
\end{cor}

In \cite[(19.7)]{F} Foxby proved an Intersection Theorem for flat
dimension.  More precisely, he showed that for $M$ an $R$-module
of finite flat dimension and of finite depth, and $R$ admitting of
a Hochster module (as is the case where $R$ is
equicharacteristic), one has:
$$\dim R \le\dim_R M+ \fd_R M.$$

Recall that the local ring $(R,\fm,k)$ is equicharacteristic if
$\mbox{char} R=\mbox{char} k$, where $\mbox{char} R$ denotes to
the characteristic of the ring $R$.  Now we extend Foxby's result
to the homological flat dimensions in the following theorem:

\begin{thm} \label{int}
Let $M$ be an $R$-module of finite depth such that $\H\fd_R
M<\infty$. Suppose that $R$ is an equicharacteristic zero ring,
then:
$$\dim R \le\dim_R M+ \H\fd_R M,$$
for $\H=\mbox{CI}$, $\mbox{G}^*$, and $\CM^*$.
\end{thm}
The proof of this theorem makes use of Lemma \ref{13} below and
the notion of the {\it Cohen-Macauley defect} ($\cmd R$) of a ring
$R$ which is defined as:
$$\cmd R:=\dim R-\depth R.$$

\begin{lem} \label{13} Let $Q\rightarrow R'$ be any $\CM^*$-deformation.
Then $\cmd R'\le \cmd Q$.
\end{lem}

\begin{proof} Suppose that $J=\ker(Q\to R')$. Since $J$ is a perfect ideal of $Q$, we
have $\pd_Q R'=\grade_Q J$. The proof of the Lemma is easily
completed by noting that $\depth Q-\depth_Q R'=\pd_Q R'=\grade_Q
J\le\hight J\le\dim Q-\dim R',$ in which the first equality
follows the Auslander-Buchsbaum formula.
\end{proof}

\noindent{\it Proof of Theorem \ref{int}.} It is sufficient to
prove the Theorem for $\H=\CM^*$. Choose a
$\CM^*$-quasi-deformation $R\rightarrow R'\leftarrow Q$ such that
$\fd_Q (M\otimes_R R')<\infty$ and $\CM^*\fd_R M=\fd_Q (M\otimes_R
R')-\fd_Q R'$. It can be seen that $Q$ is an equicharacteristic
zero ring. Since $R\rightarrow R'$ is a flat extension and
$\depth_R M<\infty$, it follows from \cite [(2.6)]{I} that
$\depth_{R'} (M\otimes_R R')<\infty$. Therefore we obtain
$\depth_Q (M\otimes_R R')<\infty$ since $Q\rightarrow R'$ is
surjective. By Lemma \ref{13} there is an inequality $\cmd
R'\le\cmd Q$. Then $\cmd R+ \cmd R'/\fm R'\le\cmd Q$. So we have:
\begin{align*}
\dim R\le & \cmd Q-\cmd R'/\fm R'+\depth R \\[1ex]
      = & \dim Q-\dim R'/\fm R'-\depth Q+\depth R+\depth R'/\fm R' \\[1ex]
      = & \dim Q-\dim R'/\fm R'-\depth Q+\depth R' \\[1ex]
      = & \dim Q-\dim R'/\fm R'-\pd_Q R' \\[1ex]
      \le &\dim_Q (M\otimes_R R')+\fd_Q(M\otimes_R R')-\fd_Q
      R'-\dim R'/\fm R' \\[1ex]
      = & \dim_{R'} (M\otimes_R R')+\CM^*\fd_R M-\dim R'/\fm R' \\[1ex]
      = & \dim_R M+ \CM^*\fd_R M,
\end{align*}
where the third equality holds by the Auslander-Buchsbaum formula;
and the second inequality holds from Foxby's Theorem
\cite[(19.7)]{F}. To prove the fifth equality assume that $M$ is
the direct union of finite submodules $M_i$ of $M$ (for $i$ in a
directed set $I$). Then
$$\dim_R M=\sup\{\dim_R M_i|i\in I\}.$$
So we get that $M\otimes_RR'$ is the direct union of
$M_i\otimes_RR'$. Consequently by the above observation we have:
\begin{align*}
\dim_{R'}(M\otimes_RR')= & \sup\{\dim_{R'}(M_i\otimes_RR')|i\in I\} \\[1ex]
                     = & \sup\{\dim_R M_i+ \dim R'/\fm R'|i\in I\}\\[1ex]
                     = & \sup\{\dim_R M_i|i\in I\}+ \dim R'/\fm R' \\[1ex]
                     = & \dim_R M+ \dim R'/\fm R',
\end{align*}
where the second equality follows from \cite [(A.11)]{BH}.
\hfill$\square$

\section{Cohen-Macaulay flat
dimension}

In this section we introduce the notion of \emph{Cohen-Macaulay flat
dimension} denoted by $\CMfd$. For a finite $R$-module $M$ it
coincides with the Cohen-Macaulay dimension $\cmdim_RM$ of Gerko.
And we show that, for an $R$-module $M$ we have the following
sequence of inequalities
$$\rfd_RM\le\CMfd_RM\le\gfd_RM\le\mbox{CI}\fd_RM\le\fd_RM,$$
with equality to the left of any finite number.

\begin{defn} Let $M\neq0$ be an $R$-module. The \textit{Cohen-Macaulay flat dimension} of $M$, is
defined as:
$$\CMfd_RM:=\inf\{\gfd_Q M' -\gfd_Q R'| \text{ }R\rightarrow R'\leftarrow Q
\text{ is a }\CM\text{-quasi-deformation}\}.$$ We complement this by
$\CMfd_R 0=-\infty$.
\end{defn}

\begin{rem} By taking the trivial $\CM$-quasi-deformation $R\rightarrow
R\leftarrow R$, one has $\CMfd_R M\le\gfd_RM$, and using Proposition
\ref{R1} we have, when $\CMfd_R M<\infty$, then $\CMfd_R M=\rfd_RM$.
\end{rem}

Notice that there is a notion of Cohen-Macaulay flat dimension in
\cite{HJ} which is different with ours. Before proceeding any
further it is necessary to investigate the effect of change of ring
on various notions of homological flat dimensions.

\begin{prop}\label{P1} Let $M$ be an $R$-module. Let $R\rightarrow R'$ be a local flat extension, and $M'=M\otimes_R
R'$. Then
$$\H\fd_R M\le\H\fd_{R'} M'$$
with equality when $\H\fd_{R'} M'$ is finite, for $\H=\mbox{CI}$,
$\mbox{G}^*$, $\CM^*$, and $\CM$.
\end{prop}

\begin{proof} We prove the result for Cohen-Macaulay flat dimension
and the proof of the other cases are similar to this one, so we omit
them. Suppose that $\CMfd_{R'} M'<\infty$, and let $R'\rightarrow
R''\leftarrow Q$ be a $\CM$-quasi-deformation with $\gfd_Q
M''<\infty$, where $M''=M'\otimes_{R'} R''$. Since $R\rightarrow R'$
and $R'\rightarrow R''$ are flat extensions, the local homomorphism
$R\rightarrow R''$ is also flat. Hence $R\rightarrow R''\leftarrow
Q$ is a $\CM$-quasi-deformation with $\gfd_Q (M\otimes_R
R'')<\infty$. It follows that $\CMfd_R M$ is finite. Now by Theorem
\ref{CH} and \cite[(8.5)]{IS}, we have
$$\CMfd_R M=\rfd_R M=\rfd_{R'} M'=\CMfd_{R'} M'.$$
\end{proof}

\begin{prop}\label{P2}
Let $\widehat{R}$ be the completion of $R$ relative to the
$\fm$-adic topology. Then
$$\H\fd_R M=\H\fd_{\widehat{R}} (M\otimes_R \widehat{R}),$$
for $\H=\mbox{CI}$, $\mbox{G}^*$, and $\CM^*$.
\end{prop}

\begin{proof} We prove the result for $\CM^*\fd$ and the proof of the other cases are similar to this one.
If $\CM^*\fd_R M=\infty$, then we obtain that
$\CM^*\fd_{\widehat{R}}(M\otimes_R \widehat{R})=\infty$ by \ref{P1}.
Now assume that $\CM^*\fd_R M<\infty$. It is sufficient to prove
that $\CM^*\fd_{\widehat{R}} (M\otimes_R \widehat{R})$ is finite.
Because in this case we have
$$\CM^*\fd_RM=\rfd_RM=\rfd_{\widehat{R}}(M\otimes_R \widehat{R})=\CM^*\fd_{\widehat{R}}(M\otimes_R \widehat{R}),$$
in which the first and the last equalities follow from Theorem
\ref{CH}, and the middle one follows from \cite[(8.5)]{IS}.

For a $\CM^*$-quasi-deformation $R\rightarrow R'\leftarrow Q$ of
$R$, we have $\widehat{R}\rightarrow \widehat{R'}\leftarrow
\widehat{Q}$ is a $\CM^*$-quasi-deformation of $\widehat{R}$ with
respect to their maximal ideal-adic completions. Now the equalities
\begin{align*}
\fd_Q (M\otimes_R R')= & \fd_{\widehat{Q}} (M\otimes_R R'\otimes_Q
\widehat{Q})=\fd_{\widehat{Q}} (M\otimes_R
\widehat{R'}) \\[1ex]
      = & \fd_{\widehat{Q}} (M\otimes_R (\widehat{R}\otimes_{\widehat{R}}
\widehat{R'}))=\fd_{\widehat{Q}} ((M\otimes_R
\widehat{R})\otimes_{\widehat{R}} \widehat{R'}),
\end{align*}
show that $\fd_{\widehat{Q}} ((M\otimes_R
\widehat{R})\otimes_{\widehat{R}} \widehat{R'})$ is finite which
imply that $\CM^*\fd_{\widehat{R}} (M\otimes_R \widehat{R})$ is
finite.
\end{proof}

One of the main result of this paper is Theorem \ref{cmi} below the
proof of which strongly makes use results of Sather-Wagstaff
\cite[Theorem F]{SW} and Esmkhani and Tousi \cite[Corollary
2.6]{ET}.

\begin{thm}\label{cmi} Let $M$ be an $R$-module.
Then there is the inequality
$$\gfd_RM\le\mbox{CI}\fd_RM$$ with equality if $\mbox{CI}\fd_RM$ is a
finite number.
\end{thm}
\begin{proof} \textbf{Step 1.} Assume that $R$ admits of a dualizing complex $D$.
We can actually assume that $\mbox{CI}\fd_RM$ is finite. So that
by \cite[Theorem F]{SW}, there exists a
$\mbox{CI}$-quasi-deformation $R\to R'\leftarrow Q$ such that $Q$
is complete, the closed fibre $R'/\fm R'$ is artinian and
Gorenstein, and $\fd_Q(M\otimes_RR')$ is finite. Therefore by
Remark \ref{Re}(a), $M\otimes_RR'$ belongs to the the Auslander
class $\mathbf{A}(Q)$. On the other hand since the kernel of $Q\to
R'$ is generated by $Q$-regular elements, using \cite[Proposition
4.3]{AF3} we deduce that it is a Gorenstein local homomorphism.
Thus thanks to \cite[Corollary (7.9)]{AF2} we see that
$M\otimes_RR'$ belongs to the the Auslander class
$\mathbf{A}(R')$. Note that $R'$ is a complete local ring, so it
admits of a dualizing complex. Hence using Remark \ref{Re}(a), we
obtain that $\gfd_{R'}(M\otimes_RR')$ is finite. Since $R'/\fm R'$
is a Gorenatein local ring, by \cite[Proposition 4.2]{AF3}, we
have $R\to R'$ is a Gorenstein local homomorphism. Therefore by
\cite[Theorem 5.1]{AF3}, the complex $D\otimes_R^\mathbf{L}R'$ is
a dualizing complex of $R'$. Consequently by \cite[Theorem
5.3]{CFH}, $\gfd_RM$ is finite. Hence the equalities
$\gfd_RM=\rfd_RM=\mbox{CI}\fd_RM$ hold.

\textbf{Step 2.} Now let $R$ be any ring. Note that by Proposition
\ref{P2} we have $$\CIfd_R M=\CIfd_{\widehat{R}} (M\otimes_R
\widehat{R}).$$ Since $\widehat{R}$ admits of a dualizing complex
by Step 1 we have $$\gfd_{\widehat{R}} (M\otimes_R
\widehat{R})\le\CIfd_{\widehat{R}} (M\otimes_R \widehat{R})$$ with
equality if $\CIfd_{\widehat{R}} (M\otimes_R \widehat{R})$ is
finite. Now assume that $\CIfd_R M$ is finite. Therefore
$\gfd_{\widehat{R}} (M\otimes_R \widehat{R})$ is finite.
Consequently by \cite[Corollary 2.6]{ET}, $\gfd_R M$ is finite and
$\gfd_R M=\gfd_{\widehat{R}} (M\otimes_R \widehat{R})$. Hence
$$\gfd_R M=\gfd_{\widehat{R}} (M\otimes_R \widehat{R})=\CIfd_{\widehat{R}} (M\otimes_R
\widehat{R})=\CIfd_R M.$$ This completes the proof.
\end{proof}

\begin{cor} Let $M$ be an $R$-module. Then there is the following sequence of
inequalities
$$\rfd_RM\le\CMfd_RM\le\gfd_RM\le\mbox{CI}\fd_RM\le\fd_RM,$$
with equality to the left of any finite number.
\end{cor}

\begin{prop}\label{P3} Let $M$ be an $R$-module. For each prime ideal $\fp \in \Spec(R)$ there is an inequality
$$\H\fd_{R_{\fp }} {M_{\fp}}\le\H\fd_R M,$$
for $\H=\mbox{CI}$, $\mbox{G}^*$, $\CM^*$, and $\CM$.
\end{prop}

\begin{proof} We prove the result for $\CMfd_{R} M$, and the proof of the other cases are similar to this one. Choose a prime
ideal $\fp \in \Spec(R)$. Assume that $\CMfd_RM<\infty$ and fix a
$\CM$-quasi-deformation $R\rightarrow R'\leftarrow Q$ such that
$\gfd_Q M'<\infty$, where $M'=M\otimes_R R'$. Since $R\rightarrow
R'$ is faithfully flat extension of rings, there is a prime ideal
$\fp '$ in $R'$ lying over $\fp$. Let $\fq$ be the inverse image of
$\fp '$ in $Q$. The map $R_{\fp}\rightarrow R'_{\fp '}$ is flat, and
$R'_{\fp '}\leftarrow Q_{\fq}$ is a $\CM$-deformation. Therefore the
diagram $R_{\fp}\rightarrow R'_{\fp '}\leftarrow Q_{\fq}$ is a
$\CM$-quasi-deformation with $\gfd_{Q_{\fq}}
(M_{\fp}\otimes_{R_{\fp}} R'_{\fp '})=\gfd_{Q_{\fq}}M'_{\fq
'}\le\gfd_QM'<\infty$. Hence $\CM\fd_{R_{\fp}} M_{\fp}<\infty$. So
we obtain
\begin{align*}
\CMfd_{R_{\fp}} M_{\fp}= & \rfd_{R_{\fp}} M_{\fp} \\[1ex]
                       \le & \rfd_R M \\[1ex]
                       = & \CMfd_R M,
\end{align*}
in which the inequality holds by \cite [(2.3)]{CFF}. Thus the
desired inequality follows.
\end{proof}

\begin{lem}\label{L1}
Let $Q$ be a local ring, and let $J\subseteq I$ be ideals of $Q$.
Set $R=Q/J$. If $J$ and $I/J$ are perfect ideals of $Q$ and $R$
respectively, then $I$ is a perfect ideal in $Q$.
\end{lem}

\begin{proof} Since $\pd_Q R< \infty$ and $\pd_R Q/I< \infty$, by \cite[(3.8)]{AB1}
there is an equality $\pd_Q Q/I=\pd_R Q/I+\pd_Q R$. By our
assumption $J$ is a perfect ideal of $Q$ hence by
\cite[(2.7)]{AF1} we have $\grade_Q Q/I=\grade_R Q/I+\grade_Q R$.
Using the perfectness of $J$ in $Q$ and $I/J$ in $R$, we see that
$I$ is a perfect ideal in $Q$.
\end{proof}

\begin{prop}\label{P4}
Let ${\bf x}=x_1,\ldots,x_n$ be a sequence of elements of $\fm$,
constituting $R$- and $M$-regular elements. Set
$\overline{R}=R/({\bf x})$ and $\overline{M}=M/({\bf x})M$. Then
there are inequalities
$$\H\fd_{\overline{R}} {\overline{M}}\le\H\fd_R M,\text{ and }$$
$$\H\fd_R {\overline{M}}\le\H\fd_R M+n,$$
with equality when, $\H\fd_R M$ is finite, for $\H=\mbox{CI}$,
$\mbox{G}^*$, and $\CM^*$.
\end{prop}

\begin{proof} Since the proof for $\H=\mbox{CI}$ and
$\mbox{G}^*$ is analogous to $\CM^*$, we only prove the
proposition for $\H=\CM^*$. It is sufficient to prove the
proposition for ${\bf x}=x$ with $x$ an $R$-regular and
$M$-regular element. We may assume that $\CM^*\fd_R M<\infty$ and
choose a $\CM^*$-quasi-deformation $R\rightarrow R'\leftarrow Q$
with $\fd_R M'<\infty$, where $M'=M\otimes_R R'$. Thus $R'=Q/J$,
where $J$ is a perfect ideal of $Q$. We construct a
$\CM^*$-quasi-deformation of $\overline{R}$. Choose $y\in Q$
mapping to $x\in R'$. Since $x$ is $R$-regular, it is also
$R'$-regular due to flatness of $R'$ as an $R$-module. Set
$I=(y)+J$ and note that $I/J=xR'$ is a perfect ideal of $R'$.
Therefore by lemma \ref{L1}, $I$ is a perfect ideal in $Q$ (for
the case $\H=\mbox{G}^*$ use \cite[(2.11)]{V}). Set
$\overline{R'}= Q/I$, and note that
$\overline{R}\rightarrow\overline{R'}$ is flat because
$R\rightarrow R'$ is flat. Thus
$\overline{R}\rightarrow\overline{R'}\leftarrow Q$ is a
$\CM^*$-quasi-deformation of $\overline{R}$.

Now we show that $\fd_Q (\overline{M}\otimes_{\overline{R}}
\overline{R'})$ and $\fd_Q (\overline{M}\otimes_R R')$ are finite.
We have the following isomorphisms
$$\overline{M}\otimes_{\overline{R}}
\overline{R'}\cong \overline{M}\otimes_R \overline{R}\otimes_R
R'\cong \overline{M}\otimes_R R'.$$ Since $x$ is $M$-regular and
$R\rightarrow R'$ is flat, the exact sequence $0\rightarrow
M\stackrel{x}{\to} M \rightarrow \overline{M}\rightarrow0$ induces
an exact sequence $0\rightarrow M'\stackrel{x}{\to} M' \rightarrow
\overline{M}\otimes_R R'\rightarrow0$. So we obtain
$\overline{M}\otimes_R R'\cong M'/xM'$ and we have $\fd_Q
(M'/xM')=\fd_Q M' +1$. Hence we get $\CM^*\fd_{\overline{R}}
\overline{M}$ and $\CM^*\fd_R \overline{M}$ are finite. Now the
equalities
$$\CM^*\fd_{\overline{R}}
\overline{M}=\rfd_{\overline{R}} \overline{M}=\rfd_R M=\CM^*\fd_R
M,$$ where the second equality follows from \cite[(3.11)]{SS},
complete the proof of the first inequality in the assertion of the
Theorem.   The equalities
$$\CM^*\fd_R \overline{M}=\rfd_R \overline{M}=\rfd_{\overline{R}} \overline{M}+1=\CM^*\fd_R M+1,$$
where the second equality follows from \cite [(3.6)]{SY2} and the
third one holds by \cite[(3.11)]{SS} complete the proof.
\end{proof}

\begin{prop}\label{P5} Let ${\bf x}=x_1,\ldots,x_n$ be a $R$-regular elements.  Set $\overline{R}=R/({\bf x})$. For an
$\overline{R}$-module $M$, then there is the inequality
$$n+\CIfd_{\overline{R}} M\le\CIfd_R M,$$
with equality when $\CIfd_R M$ is finite.
\end{prop}

\begin{proof} As usual we may assume that $\CIfd_R M<\infty$ and choose a
$\mbox{CI}$-quasi-deformation $R\rightarrow R'\leftarrow Q$ with
$\fd_Q M'<\infty$, where $M'=M\otimes_R R'$. Consider
$\overline{R}\rightarrow\overline{R}\leftarrow R$ as a
$\mbox{CI}$-quasi-deformation. One checks readily that
$\overline{R}\rightarrow R''=\overline{R}\otimes_R R'\leftarrow Q$
is a $\mbox{CI}$-quasi-deformation of $\overline{R}$. From the
equalities $M\otimes_R R'=(M\otimes_{\overline{R}}
\overline{R})\otimes_R R'=M\otimes_{\overline{R}}
(\overline{R}\otimes_R R')=M\otimes_R R''$, we obtain that $\fd_Q
(M\otimes_R R'')$ and so $\CIfd_{\overline{R}} M$ are finite. Now
the equalities
$$n+\CIfd_{\overline{R}} M=n+\rfd_{\overline{R}} M=\rfd_R M=\CIfd_R M,$$
where the second one holds by \cite [(3.6)]{SY2} complete the
proof.
\end{proof}

Let $\varphi:R\rightarrow S$ be a local homomorphism of complete
local rings. Let $N$ be a finite S-module, and let $R\rightarrow
R'\rightarrow S$ be a Cohen factorization of $\varphi$ (cf.
\cite{AFH1}). The following inequalities hold:
$$\fd_R N \leq \pd_{R'} N    \leq \fd_R N  + \mbox{edim}(R'/\fm R')$$
$$\gfd_R N \leq \gd_{R'} N \leq \gfd_R N + \mbox{edim}(R'/\fm R'),$$
where $\mbox{edim}(R'/\fm R')$ is the minimal number of generators
of the maximal ideal of $R'/\fm R'$. The first inequality is by
\cite{AFH} and the latter uses the recent characterization by
Christensen, Frankild, and Holm of certain Auslander categories in
terms of finiteness of G-dimensions (cf. \cite{CFH}, and also
\cite[Theorem 8.2]{IS}).

\begin{ques}Let $\varphi:R\rightarrow S$ be a local homomorphism of complete
local rings. Let $N$ be a finite S-module and let $R\rightarrow
R'\rightarrow S$ be a Cohen factorization of $\varphi$. The question
is whether the following inequalities hold:

$$\CIfd_R N \leq \cid_{R'}N\leq \CIfd_R N  + \mbox{edim}(R'/\fm R').$$
\end{ques}

\section{The depth formula}
The point of this section is to prove the \emph{depth formula} and
note its immediate consequences.
\begin{nota} For $R$-modules $M$ and $N$ set
$$\fd_R (M,N)=\sup\{i|\Tor_i ^R (M,N)\neq 0\}.$$ In particular, if
$\Tor_n ^R (M,N)=0$ for all $n$, then $\fd_R (M,N)=-\infty$, else
$0\le\fd_R (M,N)\le\infty$. For a finite $R$-module $M$,
$\fd_R(M,k)$ is the usual flat dimension of $M$ which is also equal
to its projective dimension $\pd_RM$. Moreover for such an $M$,
$\fd_R(M,N)$ is finite for every finitely generated $N$.
\end{nota}
\begin{thm} \label{D2} Let $M$ and $N$ be $R$-modules such that $\CIfd_R M<\infty$. If $\fd_R(M,N)<\infty$, then
$$\fd_R (M,N)\geq\depth R-\depth_R M-\depth_R N$$ with equality if
and only if $\depth_R \Tor_s ^R (M,N)=0$, for $s=\fd_R (M,N)$.
\end{thm}

\begin{proof} Since $\CIfd_R M<\infty$ there is, say a codimension $c$
$\mbox{CI}$-quasi-deformation $R\rightarrow R'\leftarrow Q$, such
that $\fd_Q M'<\infty$, where $M'=M\otimes_R R'$. By codimension
$c$ we mean that the kernel of the homomorphism $Q\to R'$ is
generated by regular elements of length $c$. Choose
$\fp\in\Spec(R')$ such that it is a minimal prime ideal containing
$\fm R'$.  Thus $\fm=\fp\cap R$ and $\fp=\fq/(\bf x)$ for some
$\fq\in\Spec(Q)$, where $({\bf x})=\ker(Q\rightarrow R')$. Now the
diagram $R\rightarrow R'_{\fp}\leftarrow Q_{\fq}$ is a
$\mbox{CI}$-quasi-deformation of the same codimension as
$R\rightarrow R'\leftarrow Q$. It is clear that $\pd_Q
R'=\pd_{Q_{\fq}} R'_{\fp}$. Also we have
$$\fd_{Q_{\fq}} (M\otimes_R R'_{\fp})=\fd_{Q_{\fq}} (M\otimes_R (R'\otimes_Q Q_{\fq}))
=\fd_{Q_{\fq}} ((M\otimes_R R')\otimes_Q Q_{\fq})\le \fd_Q
M'<\infty.$$ Hence $\CIfd_R M\le\fd_{Q_{\fq}} (M\otimes_R
R'_{\fp})-\fd_{Q_{\fq}} R'_{\fp}$. Therefore we showed that
complete intersection flat dimension can be computed from
$\mbox{CI}$-quasi-deformations $R\rightarrow R'\leftarrow Q$ such
that the closed fiber $R'/\fm R'$ is artinian.

Due to faithful flatness of $R'$ we have the following equalities
in which $N'=N\otimes_R R'$
$$s=\fd_R (M,N)=\fd_{R'} (M',N').$$ Assume that $c=1$. Consider the change of rings spectral
sequence
$$\Tor_p ^{R'} (M',\Tor_q ^Q (R',N'))\Rightarrow\Tor_{p+q} ^Q
(M',N').$$ If $q>1$, then $\Tor_q ^Q (R',N')=0$ and for $q\le 1$
$\Tor_q ^Q (R',N')=N'$. Now the above spectral sequence generates
the following long exact sequence
$$\cdots\rightarrow \Tor_{i+1} ^{R'} (M',N')\rightarrow\Tor_{i-1} ^{R'} (M',N')
\rightarrow\Tor_i ^Q (M',N')\rightarrow\Tor_i ^{R'}
(M',N')\rightarrow\cdots.$$ Therefore $\Tor_{s+1} ^Q (M',N')=\Tor_s
^{R'} (M',N')$. Iterating in the same manner we have
$$\Tor_s ^{R'} (M',N')=\Tor_{s+c} ^Q (M',N').$$
So $\sup\{i|\Tor_i ^Q (M',N')\neq 0\}=s+c$. Since $\depth(R'/\fm
R')=0$ and $Q\rightarrow R'$ is surjective, the following equalities
hold:
$$\depth_Q \Tor_s ^{R'} (M',N')=\depth_{R'} \Tor_s ^{R'} (M',N')=\depth_R \Tor_s ^R (M,N),$$
and they are equal to $\depth_Q \Tor_{s+c} ^Q (M',N')$. Since $\fd_Q
M'<\infty$ it follows from \cite [(2.3)]{SY2} that
\begin{align*}
s+c\geq & \depth Q-\depth_Q M'-\depth_Q N' \\[1ex]
      = & \depth R+c-\depth_R M-\depth_R N,
\end{align*}
and equality holds if and only if $\depth_Q \Tor_{s+c} ^Q
(M',N')=0.$ Thus $$s=\fd_R (M,N)\geq\depth R-\depth_R M-\depth_R
N,$$ with equality if and only if $\depth_R \Tor_s ^R (M,N)=0.$
\end{proof}

\begin{defn} We say that $M$ and $N$ satisfy the \emph{dependency formula} over R, if
$$\fd_R (M,N)=\sup\{\depth R_{\fp}-\depth_{R_{\fp}} M_{\fp}-\depth_{R_{\fp}}
N_{\fp}|\fp\in\Supp(M)\cap\Supp(N)\}.$$
\end{defn}

\begin{cor}\label{l} Let $M$ and $N$ be $R$-modules such that $\CIfd_R M<\infty$.
If $\fd_R(M,N)<\infty$, then $M$ and $N$ satisfy the dependency
formula.
\end{cor}

\begin{proof} It is easy to see that:
$$\depth
R_{\fp}-\depth_{R_{\fp}}M_{\fp}-\depth_{R_{\fp}}N_{\fp}\le\fd_R(M,N).$$
Using Theorem \ref{D2} we have $\fd_R (M,N)=\depth
R_{\fp}-\depth_{R_{\fp}} M_{\fp}-\depth_{R_{\fp}} N_{\fp}$ if and
only if
$$\depth_{R_{\fp}}\Tor_s ^{R_{\fp}}(M_{\fp},N_{\fp})=0,$$
or equivalently if and only if $\fp\in\Ass(\Tor_s ^R (M,N))$ for
$s=\fd_R (M,N)<\infty$.
\end{proof}

In Theorem \ref{cmi} we proved that $\gfd_RM\le\CIfd_RM$ for any
$R$-module $M$. So it is natural to look for a dependency formula
for Gorenstein flat dimension. In the following proposition we
prove a dependency formula for Gorenstein flat dimension:
\begin{prop} \label{D1}
Let $M$ and $N$ be $R$-modules, such that $\gfd_R M<\infty$ and
$\id_RN<\infty$. Then then $M$ and $N$ satisfy the dependency
formula.
\end{prop}
\begin{proof}
It is clear that $\gfd_{R_{\fp}}M_{\fp}$ and $\id_{R_{\fp}}M_{\fp}$
are finite numbers, and we have
$\fd_{R_{\fp}}(M_{\fp},N_{\fp})\le\fd_R(M,N)$ for
$\fp\in\Supp(M)\cap\Supp(N)$. Now from \cite[(12.26)]{F} and \cite
[(4.4)(a)]{CH} we get:
$$\depth
R_{\fp}-\depth_{R_{\fp}}M_{\fp}-\depth_{R_{\fp}}N_{\fp}\le\fd_R(M,N),$$
with equality when $\fp\in\Ass(\Tor_s^R(M,N))$, for $s=\fd_R(M,N)$.
\end{proof}

The following corollary due to Iyengar and Sather-Wagstaff
\cite[Theorem (8.7)]{IS} is an immediate consequence of the above
proposition.

\begin{cor} Let $M$ be an $R$-module such that
$\gfd_RM<\infty$. Then
$$\sup\{i|\Tor_i
^R (M,E_R(k))\neq 0\}=\depth R-\depth_RM,$$ where $E_R(k)$ denotes
to the injective envelope of $k$.
\end{cor}
\begin{cor} \label{D3}
Let $(R,\fm)$ be a complete local ring and $N$ and $M$ be
$R$-modules with $M$ finite.
\begin{itemize}
\item[(a)] If $\gfd_RN<\infty$, and $fd_RM<\infty$, then
$$\sup\{i|\Ext^i_R(N,M)\neq0\}=\depth R-\depth_RN.$$

\item[(b)] If $CIfd_RN<\infty$, then
$$\sup\{i|\Ext^i_R(N,M)\neq0\}=\depth R-\depth_RN,$$
provided that the left hand side is a finite number.

\end{itemize}
\end{cor}
\begin{proof}
(a) Set $E = E_R(k)$, the injective envelope of $k$. Since R is
complete we have $R =\Hom_R(E,E)$. Therefore we have
\begin{align*}
\Ext^i_R(N,M)\cong & \Ext^i_R(N,M\otimes_RR)\cong\Ext^i_R(N,M\otimes_R\Hom_R(E,E)) \\[1ex]
      \cong & \Ext^i_R(N,\Hom_R(\Hom_R(M,E),E)) \\[1ex]
      \cong & \Hom_R(\Tor^R_i(N,\Hom_R(M,E),E).
\end{align*}
Consequently we have
$$\sup\{i|\Ext^i_R(N,M)\neq0\}=\fd_R(N,\Hom_R(M,E)).$$
Since $\fd_RM<\infty$ we have $\id_R \Hom_R(M,E)<\infty$. Now
Proposition \ref{D1}, gives the result.

(b) Similarly to that of part (a) one has
$$\sup\{i|\Ext^i_R(N,M)\neq0\}=\fd_R(N,\Hom_R(M,E)),$$
which is equal to $\depth R-\depth_RN$ by Corollary \ref{l} and
\cite [Lemma 2.2]{Y1}.
\end{proof}

The following example shows that the completeness assumption of $R$
is crucial in part (a) of the above corollary.

\begin{exam}
Let $(R,\fm)$ be a local domain which is not complete with respect
to the $\fm$-adic topology. In \cite [(3.3)]{AER} it is shown that
$\Hom_R(\widehat{R},R)=0$. Therefore, when $N=\widehat{R}$ and
$M=R$.  It is clear that the right hand side of the first equality
equal to zero which is not equal to the left hand side.
\end{exam}

\begin{cor}
Let $M$ and $N$ be $R$-modules;
\begin{itemize}
\item[(a)] If $\CIfd_RM<\infty$ then the following are equivalent:
\begin{itemize}
\item[(i)] $\Tor^R_n(N,M)=0$ $n\gg0$.

\item[(ii)] $\Tor^R_n(N,M)=0$ $n>\CIfd_RM$.

\end{itemize}

\item[(b)] If $R$ is a complete local ring and $\CIfd_RN<\infty$, and $M$ a finite $R$-module, then the
following are equivalent:

\begin{itemize}

\item[(i)] $\Ext_R^n(N,M)=0$ $n\gg0$.

\item[(ii)] $\Ext_R^n(N,M)=0$ $n>\depth R-\depth_RM$.

\end{itemize}
\end{itemize}
\end{cor}
\begin{proof}
(a) If for all integer $n$, $\Tor^R_n(M,N)=0$, then the assertion
holds. So assume for some integer $\ell$, $\Tor^R_{\ell}(M,N)\neq0$.
Therefore $s=\fd_R(M,N)<\infty$. Now by Theorem \ref{D2}, $s=\depth
R_{\fp}-\depth_{R_{\fp}}M_{\fp}-\depth_{R_{\fp}}N_{\fp}$ for some
$\fp\in\Supp(M)\cap\Supp(N)$. Now choose an integer
$n>\CIfd_RM=\rfd_RM\geq\depth R_{\fp}-\depth_{R_{\fp}}M_{\fp}\geq
s$. Therefore $\Tor^R_n(N,M)=0$.

(b) It follows easily from Corollary \ref{D3}(b).
\end{proof}

\section{The finite case}

In this section $M$ is a finite $R$-module. We study the behavior
of the new homological flat dimensions and especially of the
(upper) Cohen-Macaulay flat dimension. Since for a finite
$R$-module $M$, there is the equality $\fd_RM=\pd_RM$ and
$\gfd_RM=\gd_RM$ by \cite{EJT}, we have
$$\CM^*\fd_R M=\inf\{\pd_Q M' -\pd_Q R'| \text{ }R\rightarrow R'\leftarrow
Q \text{ is a }\CM^*\text{-quasi-deformation}\},$$ $\CMfd_R
M=\cmdim_RM$, $\mbox{G}^*\fd _R M=\mbox{G}^*\text{-}\dim_RM$, and
$\mbox{CI}\fd_RM=\cid_RM$.

\begin{rem} \label{C2} It can be seen that if $\CM^*\fd_R M<\infty$, then there is the
equality $$\CM^*\fd_R M+\depth_R M=\depth R.$$
\end{rem}

Let $\mbox{Syz}_n ^R (M)$ to denote the $n$-th syzygy module of
$M$. Then by an argument similar to that of \cite[(2.5)]{V} we
have the following proposition:

\begin{prop} For each $n\geq 0$ there is the
 equality
 $$\CM^*\fd_R \mbox{Syz}_n ^R (M)=\max\{\CM^*\fd_R M
 -n,0\}.$$
\end{prop}

\begin{thm} \label{CHAR} The following conditions are equivalent:
\begin{itemize}
\item[(i)] The ring $R$ is Cohen-Macaulay.

\item[(ii)] $\CM^*\fd_R M<\infty$ for every not necessarily finite $R$-module $M$.

\item[(iii)] $\CM^*\fd_R M<\infty$ for every finite $R$-module $M$.

\item[(iv)] $\CM^*\fd_R M=0$ for every finite $R$-module $M$ with $\depth_R M\geqslant \depth R$.

\item[(v)] $\CM^*\fd_R k=\depth R$.

\item[(vi)] $\CM^*\fd_R k<\infty$.

\end{itemize}

\end{thm}

\begin{proof} $(i)\Rightarrow (ii)$ Let $\widehat{R}$ be
the $\fm$-adic completion of $R$. Since $R$ is Cohen-Macaulay, so
is $\widehat{R}$. Therefore by Cohen's structure theorem,
$\widehat{R}$ is isomorphic to $Q/J$, where $Q$ is a regular local
ring. By Cohen-Macaulay-ness of $\widehat{R}$ and regularity of
$Q$, the ideal $J$ is perfect.  Thus $R\rightarrow
\widehat{R}\leftarrow Q$ is a $\CM^*$-quasi-deformation. Since $Q$
is regular $\fd_Q(M\otimes_R\widehat{R})$ is finite. Thus
$\CM^*\fd_R M$ is finite.

 $(ii)\Rightarrow (iii)$ is trivial.

 $(iii)\Rightarrow (iv)$ follows by applying Remark \ref{C2}
to the $R$-module $M$.

 $(iv)\Rightarrow (v)$ by \cite [(1.3.7)]{BH} we have
$\depth_R \mbox{Syz}_n ^R (k)\geq\min(n, \depth R).$ In particular,
if we choose $n\geq\depth R$ we get $\CM^*\fd_R M=0$. Thus by Remark
\ref{C2} $\CM^*\fd_R k=\depth R$.

 $(v)\Rightarrow (vi)$ is trivial.

 $(vi)\Rightarrow (i)$ follows from $\cmdim_Rk\le\CM^*\fd_Rk$ and \cite[Theorem (3.9)]{G}.
\end{proof}
One can actually state similar theorems for upper Gorenstein flat
and complete intersection flat dimensions.

As a consequence of the New Intersection Theorem of Peskine and
Szpiro \cite{PS}, Hochster \cite{HH} and P. Roberts \cite{R1} and
\cite{R2} we have:
$$(*)\qquad\cmd R\le\cmd_RM.$$

The New Intersection Theorem is not true for $\mbox{CI}$-dimension,
$\mbox{G}^*$-dimension, $\mbox{G}$-dimension,  and $\CM$-dimension,
see Examples \cite [(3.2)]{SY1} and \cite [(2.20)]{Y}. But the
inequality $(*)$ holds for $\mbox{G}^*$-dimension, see
\cite[(2.1)]{SY1}. For $\mbox{G}$-dimension and
$\mbox{CM}$-dimension we do not know whether the inequality $(*)$
holds.  However it holds for upper Cohen-Macaulay flat dimension as
shown by the following theorem:

\begin{thm} Let $M$ be a finite $R$-module with finite
upper Cohen-Macaulay dimension. Then $\cmd R\le\cmd_R M$.
\end{thm}

\begin{proof} Since $\CM^*\fd_R M<\infty$, there exists a
$\CM^*$-quasi-deformation $R\rightarrow R'\leftarrow Q$ such that
$\pd_Q M'<\infty$, where $M'=M\otimes_R R'$. From the surjectivity
of $Q\to R'$ we have $\cmd_Q M'=\cmd_{R'} M'$. Because $R\rightarrow
R'$ is a flat extension, the following (in)equalities hold:
\begin{align*}
\cmd R+\cmd_R R'/\fm R'= & \cmd R'\le\cmd Q\le\cmd_Q M' \\[1ex]
                       = & \cmd_R M+\cmd R'/\fm R',
\end{align*}
where the first inequality holds by Lemma \ref{13}, and the second
one is by the New Intersection Theorem. This gives us the desired
inequality.
\end{proof}

\begin{cor} If $M$ is a Cohen-Macaulay module with $\CM^*\fd_R M<\infty$, then the
base ring $R$ is Cohen-Macaulay.
\end{cor}

\section{Homological Injective Dimensions}
It is well known that flat dimension and injective dimension are
dual of each other. In particular there are the following
equalities:
$$\fd_RM^\vee=\id_RM\text{ and }\id_RM^\vee=\fd_RM,$$ where
$M^\vee=\Hom_R(M,E_R(k))$ and $E_R(k)$ is the injective envelope of
$k$ over $R$. In this section we introduce dual of the complete
intersection flat dimension and the Cohen-Macaulay flat dimension.
\begin{defn} Let $M\neq0$ be an $R$-module. The \textit{complete
intersection injective dimension}, \textit{upper Gorenstein
injective dimension}, and \textit{upper Cohen-Macaulay injective
dimension} of $M$, are defined as:
$$\mbox{CI}\id_RM:=\inf\{\id_Q M' -\fd_Q R'| \text{ }R\rightarrow R'\leftarrow
Q \text{ is a }\mbox{CI}\text{-quasi-deformation}\}$$
$$\mbox{G}^*\id_RM:=\inf\{\id_Q M' -\fd_Q R'| \text{ }R\rightarrow R'\leftarrow
Q \text{ is a }\mbox{G}^*\text{-quasi-deformation}\}$$
$$\CM^*\id_RM:=\inf\{\id_Q M' -\fd_Q R'| \text{ }R\rightarrow R'\leftarrow Q
\text{ is a }\CM^*\text{-quasi-deformation}\},$$ respectively.
We complement this by $\H\fd_R 0=-\infty$ for $\H=\mbox{CI}$,
$\mbox{G}^*$, and $\CM^*$.
\end{defn}

In \cite[Definition 2.6]{SW} Sather-Wagstaff introduced the upper
complete intersection injective dimension of an $R$-module $M$ as
$$\mbox{CI}^*\id_RM:=\inf\left\{\id_Q M' -\fd_Q R'\bigg| \begin{array}{l} R\rightarrow R'\leftarrow Q\text{ is a
}\mbox{CI}\text{-quasi-deformation }\\\text{ such that }R'\text{ has
Gorenstein formal }\\\text{ fibre and }R'/\fm R'\text{ is Gorenstein
} \end{array} \right\}.$$

We use the upper complete intersection injective dimension for a
dual version of Theorem \ref{cmi}.  The following theorem shows
that the upper Cohen-Macaulay injective dimension characterizes
Cohen-Macaulay local rings.

\begin{thm}\label{OH} The following conditions are equivalent.
\begin{itemize}
\item[(i)] The ring $R$ is Cohen-Macaulay.

\item[(ii)] $\CM^*\id_R M<\infty$ for every $R$-module $M$.

\item[(iii)] $\CM^*\id_R M<\infty$ for every finite $R$-module $M$.

\item[(iv)] $\CM^*\id_R k<\infty$.

\end{itemize}

\end{thm}

\begin{proof} $(i)\Rightarrow (ii)$ Let $\widehat{R}$ be
the $\fm$-adic completion of $R$. Since $R$ is Cohen-Macaulay, so
is $\widehat{R}$. Therefore by Cohen's structure theorem,
$\widehat{R}$ is isomorphic to $Q/J$, where $Q$ is a local regular
ring. Hence due to Cohen-Macaulay-ness of $\widehat{R}$ and
regularity of $Q$, the ideal $J$ is perfect.  Thus $R\rightarrow
\widehat{R}\leftarrow Q$ is a $\CM^*$-quasi-deformation. Since $Q$
is regular $\id_Q(M\otimes_R\widehat{R})$ is finite, so
$\CM^*\id_R M$ is finite.

 $(ii)\Rightarrow (iii)$ and $(iii)\Rightarrow (iv)$ are trivial.

 $(iv)\Rightarrow (i)$ Suppose $\CM^*\id_R k<\infty$. So that there
 exists a $\CM^*$-quasi-deformation $R\rightarrow R'\leftarrow Q$,
 such that $\id_Q(k\otimes_RR')$ is finite. It is clear that $k\otimes_RR'$
 is a finite $Q$-module. Consequently $Q$ is a
 Cohen-Macaulay ring by the Bass Theorem. We plan to show that $R'$ is a
 Cohen-Macaulay ring. Let $I=\ker(Q\rightarrow R')$ which is perfect by definition.
 We have
\begin{align*}
\hight I= &\grade(I,Q) \\[1ex]
        = & \pd_QR' \\[1ex]
        = & \depth Q-\depth_QR' \\[1ex]
        = & \depth Q-\depth R' \\[1ex]
        = & \dim Q-\depth R' \\[1ex]
        = & \hight I+\dim R'-\depth R',
\end{align*}
in which the equalities follow from Cohen-Macaulay-ness of $Q$;
perfectness of $I$; Auslander-Buchsbaum formula;
\cite[(1.2.26)]{BH}; Cohen-Macaulay-ness of $Q$; and \cite[Page
250]{M} respectively. Therefore we obtain that $\dim R'-\depth
R'=0$, that is $R'$ is Cohen-Macaulay. Now \cite[Theorem
(2.1.7)]{BH} gives us the desired result.
\end{proof}

In the same way one can show that the upper Gorenstein injective
dimension detects the Gorenstein property and the complete
intersection injective dimension detects the complete intersection
propoerty of local rings.

The proof of the above theorem says some thing more, viz., a local
ring $R$ is Cohen-Macaulay if and only if there exists a finite
$R$-module of finite upper Cohen-Macaulay injective dimension. In
other words every finite $R$-module is a test module for the
Cohen-Macaulay property of a local ring.  So we state the
following corollary, which is analogous to the definition of a
Gorenstein ring:

\begin{cor} A local ring $R$ is Cohen-Macaulay if and only if
$\CM^*\id_RR<\infty$.
\end{cor}

In \cite[Theorem (4.5)]{FF} Foxby and Frankild proved that if a ring
admits a cyclic module of finite Gorenstein injective dimension,
then the base ring is Gorenstein. Parallel to their result we have
the the following proposition.

\begin{prop} If $\mbox{G}^*id_RC<\infty$ for a cyclic $R$-module $C$, then $R$ is a Gorenstein
local ring.
\end{prop}

\begin{proof} There is a $\mbox{G}^*$-quasi-deformation
$R\rightarrow R'\leftarrow Q$ such that $\id_Q(C\otimes_RR')$ is
finite. Since $C\otimes_RR'$ is a cyclic $R'$-module and $R'$ is a
cyclic module over $Q$, we see that $C\otimes_RR'$ is a cyclic
module over $Q$. So that $Q$ is Gorenstein by \cite{PZ}. Hence $R'$
is a Gorenstein ring because the kernel of $Q\rightarrow R'$ is a
Gorenstein ideal. Consequently $R$ is Gorenstein.
\end{proof}
\begin{lem} \label{art} There is an equality
$$\H\id_RM=\inf\left\{\id_Q M' -\fd_Q R'\bigg| \begin{array}{l} R\rightarrow R'\leftarrow Q\text{ is an
}\H\text{-quasi-deformation such }
\\\text{ that the closed fibre of }R\to R'\text{ is artinian } \end{array} \right\},$$
for $\H=\mbox{CI}^*$, $\mbox{CI}$, $\mbox{G}^*$, and $\CM^*$.
\end{lem}

\begin{proof} We prove the lemma for $\H=\CM^*$ only since the other cases are similar.
Let $R\rightarrow R'\leftarrow Q$ be an $\CM^*$-quasi-deformation.
Choose $\fp\in\Spec(R')$ such that it is a minimal prime ideal
containing $\fm R'$; thus $\fm=\fp\cap R$ and $\fp=\fq/J$ for some
$\fq\in\Spec(Q)$, where $J=\ker(Q\rightarrow R')$. Now the diagram
$R\rightarrow R'_{\fp}\leftarrow Q_{\fq}$ is a
$\CM^*$-quasi-deformation.  It is clear that $\pd_Q
R'=\pd_{Q_{\fq}} R'_{\fp}$.  Also we have
$$\id_{Q_{\fq}} (M\otimes_R R'_{\fp})=\id_{Q_{\fq}} (M\otimes_R (R'\otimes_Q Q_{\fq}))
=\id_{Q_{\fq}} (M'\otimes_Q Q_{\fq})\le\id_Q M'<\infty.$$ Hence
$\CM^*\id_R M\le\id_{Q_{\fq}} (M\otimes_R R'_{\fp})-\pd_{Q_{\fq}}
R'_{\fp}$. So the proof in complete.
\end{proof}

Recall that $\width_RM=\inf\{i|\Tor_i^R(M,k)\neq 0\}.$ It is the
dual notion for $\depth_RM$. In particular by \cite[(4.8)]{CFF} we
have $\width_RM= \depth_R\Hom_R(M,E_R(k))$, where $E_R(k)$ denote
for injective envelope of $k$ over $R$.

The \emph{Chouinard injective dimension} is denoted by
$\mbox{Ch}\id_RM$ and is defined as,
$$\mbox{Ch}\id_RM:=\sup \{\depth R_\fp - \width _{R_{\fp}} M_\fp | \fp\in \Spec
(R)\}.$$

\begin{rem} We do not introduce Chouinard flat dimension because it
coincides with $\rfd$ since by \cite[(2.4)(b)]{CFF}, we have
$$\rfd_RM:=\sup \{\depth R_\fp - \depth _{R_{\fp}} M_\fp |
\fp\in \Spec (R)\}.$$
\end{rem}

It is proved in \cite{C} that for an $R$-module $M$,
$\mbox{Ch}\id_RM$ is a refinement of $\id_RM$, that is
$$\mbox{Ch}\id_RM\le\id_RM,$$ with equality if $\id_RM$ is finite.

In the Theorem \ref{T} we partly extend this relation for our
homological injective dimensions. As of the writing of this paper,
the authors do not know if the equality holds in general. Before
stating the theorem we need some lemmas.

\begin{lem} \label{CH2} Suppose that $Q\rightarrow S$ is a surjective local homomorphism and $N$ is an
$S$-module. Then we have
$$\width_SN=\width_QN.$$
\end{lem}
\begin{proof} We have the following equalities:
\begin{align*}
\width_SN= & \depth_S\Hom_S(N,E_S(k)) \\[1ex]
         = & \depth_S\Hom_S(N,\Hom_Q(S,E_Q(k))) \\[1ex]
         = & \depth_S\Hom_Q(N,E_Q(k)) \\[1ex]
         = & \depth_Q\Hom_Q(N,E_Q(k)) \\[1ex]
         = & \width_QN,
\end{align*}
where the first one is by \cite[(4.8)]{CFF}; the second one is by
\cite[(10.1.15)]{BS}; the third one is by adjointness of $\Hom$
and tensor; the fourth one is true since $Q\rightarrow S$ is
surjective; while the last one is again by \cite[(4.8)]{CFF}. Here
we used $k$ for the residue fields of $Q$ and $S$, and $E_Q(k)$
and $E_S(k)$ for the injective envelopes of $k$ over respectively
$Q$ and $S$.
\end{proof}

Dualizing the proof of Proposition \ref{R1} and using the above
lemma one easily shows

\begin{prop} \label{CH5} Let $Q\rightarrow S$ be a $\CM$-deformation, and $N$ be an
$S$-module. Then there is the equality:
$$\mbox{Ch}\id_SN+\gd_Q S= \mbox{Ch}\id_QN.$$
\end{prop}

\begin{lem} \label{CH41} Suppose that $(R,\fm,k)\rightarrow (S,\fn,l)$ is a local ring homomorphism, and $M$ is an
$R$-module. Then we have
$$\width_S(M\otimes_RS)=\width_RM.$$
\end{lem}
\begin{proof} Let $F_M$ be a flat resolution of $M$ over $R$. Therefore
$F_M\otimes_RS$ is a flat resolution of $M\otimes_RS$ over $S$. So
we have
\begin{align*}
\width_S(M\otimes_RS)= & \inf\{i|\Tor_i^S(M\otimes_RS,l)\neq0\} \\[1ex]
                     = & \inf\{i|\H_i((F_M\otimes_RS)\otimes_Sl)\neq0\} \\[1ex]
                     = & \inf\{i|\H_i(F_M\otimes_Rl)\neq0\} \\[1ex]
                     = & \inf\{i|\Tor_i^R(M,l)\neq0\} \\[1ex]
                     = & \inf\{i|\Tor_i^R(M,k)\neq0\} \\[1ex]
                     = & \width_RM.
\end{align*}
\end{proof}

\begin{lem} \label{CH4} Let $R\rightarrow S$ be a flat local homomorphism and let $M$ be an
$R$-module. Then
$$Ch\id_RM\le Ch\id_S(M\otimes_RS).$$
\end{lem}

\begin{proof} Let $\fp\in\Spec(R)$ such that $\mbox{Ch}\id_RM= \depth R_{\fp}-\width_{R_{\fp}}M_{\fp}$. Let
$\fq\in\Spec(S)$ contain $\fp S$ minimally. Since $R\rightarrow S$
is a flat local homomorphism we have $\fp=\fq\cap R$ and
$\hight\fp=\hight\fq$. Hence:
\begin{align*}
\mbox{Ch}\id_RM= & \depth R_{\fp}-\width_{R_{\fp}}M_{\fp} \\[1ex]
               = & \depth S_{\fq}-\width_{S_{\fq}}(M_{\fp}\otimes_{R_{\fp}}S_{\fq}) \\[1ex]
               = & \depth S_{\fq}-\width_{S_{\fq}}(M\otimes_RS)_{\fq} \\[1ex]
               \le & \mbox{Ch}\id_S(M\otimes_RS),
\end{align*}
in which the second equality holds by Lemma \ref{CH41} and the fact
that $R_{\fp}\rightarrow S_{\fq}$ has artinian closed fibre.
\end{proof}

\begin{thm} \label{T} Suppose that $M$ is an $R$-module such that $H\id_RM<\infty$ for $\H=\mbox{CI}^*$,
$\mbox{CI}$, $\mbox{G}^*$, or $\CM^*$. Then there is the
inequality
$$Ch\id_RM\le\H\id_RM,$$ and if $M$ is a finite module we have
$$Ch\id_RM=\H\id_RM=\depth R.$$
\end{thm}

\begin{proof} We prove the theorem for $\H=\CM^*$ and the other cases are similar.
Choose by Lemma \ref{art} a $\CM^*$-quasi-deformation
$R\rightarrow R'\leftarrow Q$, such that $\CM^*\id_RM=\id_Q
M'-\fd_QR'$, where $M'=M\otimes_R R'$, and the closed fibre of
$R\rightarrow R'$ is artinian. Hence we have
\begin{align*}
\CM^*\id_RM= & \id_Q M'-\fd_QR' \\[1ex]
        = & \mbox{Ch}\id_QM'-\fd_QR' \\[1ex]
        = & \mbox{Ch}\id_{R'}M'\geq \mbox{Ch}\id_RM,
\end{align*}
in which the second equality comes by \cite{C}, and the third one by
Proposition \ref{CH5}; while the inequality is by Lemma \ref{CH4}.

Now let $M$ be a finite $R$-module, therefore $M'$ is a finite
$Q$-module. So by the Bass Theorem \cite[(18.9)]{M}, and the
Auslander-Buchsbaum formula, and the fact that the closed fibre of
$R\rightarrow R'$ is artinian we have:
$$\CM^*\id_RM=\id_Q
M'-\fd_QR'=\depth Q-\depth Q+\depth R'=\depth R'=\depth R.$$ On
the other hand, since $\CM^*\id_RM<\infty$ and $M$ is finite, $R$
is a Cohen-Macaulay ring. Thus using \cite[(3.6)]{SS} we see that
$\mbox{Ch}\id_RM=\depth R$. Hence
$\CM^*\id_RM=\mbox{Ch}\id_RM=\depth R$.
\end{proof}

\begin{cor} Let $M$ be an $R$-module. Then we have the following chain of inequalities:
$$\mbox{Ch}\id_RM\le\CM^*\id_R M\le\mbox{G}^*id_R M\le\CIid_RM\le\mbox{CI}^*\id_RM\le\id_RM,$$
with equality to the left of any finite number for finite modules
or, if $\id_RM<\infty$ for arbitrary module $M$.
\end{cor}
Now we prove the dual result of Theorem \ref{cmi}.

\begin{thm}\label{W} Suppose that $R$ has a dualizing complex $D$, and $M$ is an $R$-module.
Then there is an inequality
$$\gid_RM\le\mbox{CI}^*\id_RM.$$
\end{thm}
\begin{proof} We can actually assume that $\mbox{CI}^*\id_RM$ is finite.
So that by \cite[Proposition 3.5]{SW}, there exists a
$\mbox{CI}$-quasi-deformation $R\to R'\leftarrow Q$ such that $Q$
is complete, the closed fibre $R'/\fm R'$ is artinian and
Gorenstein and $\id_Q(M\otimes_RR')$ is finite. Therefore by
Remark \ref{Re}(b), $M\otimes_RR'$ belongs to the the Bass class
$\mathbf{B}(Q)$. By the same argument as in the proof of Theorem
\ref{cmi} one can show that $M\otimes_RR'$ belongs to the Bass
class $\mathbf{B}(R')$. Consequently Remark \ref{Re}(b), gives us
that $\gid_{R'}(M\otimes_RR')<\infty$. Since $R'/\fm R'$ is a
Gorenatein local ring, by \cite[Proposition 4.2]{AF3}, we have
$R\to R'$ is a Gorenstein local homomorphism. Therefore by
\cite[Theorem 5.1]{AF3}, the complex $D\otimes_R^\mathbf{L}R'$ is
a dualizing complex of $R'$. Consequently \cite[Theorem 5.3]{CFH}
gives finiteness of $\gid_RM$. Hence using \cite[Theorem 6.8]{CFH}
we have $\gid_RM=\mbox{Ch}\id_RM\le\mbox{CI}^*\id_RM$ as desired.
\end{proof}

\noindent By the above theorem we have
$$\mbox{Ch}\id_RM\le\gid_R M\le\mbox{CI}^*\id_RM\le\id_RM.$$ Now we
define a Cohen-Macaulay injective dimension to complete this
sequence of inequalities. Notice that there is a notion of
Cohen-Macaulay injective dimension in \cite{HJ} which is different
with ours.

\begin{defn} Let $M\neq0$ be an $R$-module. The \textit{Cohen-Macaulay injective dimension} of $M$, is
defined by ,
$$\CMid_RM:=\inf\{\gid_Q M' -\gfd_Q R'| \text{ }R\rightarrow R'\leftarrow Q
\text{ is a }\CM\text{-quasi-deformation}\}.$$ We complement this by
$\CMid_R 0=-\infty$.
\end{defn}

Therefore by taking the trivial $\CM$-quasi-deformation
$R\rightarrow R\leftarrow R$, one has $\CMid_R M\le\gid_RM$.
Suppose that $R$ has a dualizing complex. Then by \cite[Proposion
5.5]{CFH} one has $\gid_{R_{\fp}}M_{\fp}\le\gid_RM$ for each prime
ideal $\fp$ of $R$. So that by the proof of Lemma \ref{art}, one
has
$$\CMid_RM=\inf\left\{\gid_Q M' -\gfd_Q R'\bigg| \begin{array}{l} R\rightarrow R'\leftarrow Q\text{ is a
}\CM\text{-quasi-deformation }
\\\text{ such that the closed fibre of }\\R\to R'\text{ is artinian } \end{array} \right\}.$$
Therefore as in Theorem \ref{T} one can show that
$\mbox{Ch}\id_RM\le\CM\id_RM$. Hence when the ring $R$ admits a
dualizing complex, then there is the following sequence of
inequalities
$$\mbox{Ch}\id_RM\le\CMid_RM\le\gid_RM\le\mbox{CI}^*\id_RM\le\id_RM,$$
with equality to the left of any finite number for finite modules
or, if $\id_RM$ or $\gid_RM$ if finite, for arbitrary module $M$.
\begin{thm} The following conditions are equivalent.
\begin{itemize}
\item[(i)] The ring $R$ is Cohen-Macaulay.

\item[(ii)] $\CMid_R M<\infty$ for every $R$-module $M$.

\item[(iii)] $\CMid_R M<\infty$ for every finite $R$-module $M$.

\item[(iv)] $\CMid_R k<\infty$.

\end{itemize}

\end{thm}

\begin{proof} $(i)\Rightarrow (ii)$ It follows by the inequality
$\CMid_RM\le\CM^*\id_RM$ and Theorem \ref{OH}.

 $(ii)\Rightarrow (iii)$ and $(iii)\Rightarrow (iv)$ are trivial.

 $(iv)\Rightarrow (i)$ Suppose $\CM^*\id_R k<\infty$. So there
 exists a $\CM$-quasi-deformation $R\rightarrow R'\leftarrow Q$,
 such that $\gid_Q(k\otimes_RR')$ is finite. Since $k\otimes_RR'$ is a
 cyclic module over $Q$, we see that $Q$ is a Gorenstein ring by \cite[Theorem 4.5]{FF}.
 The rest of proof is the same as that of Theorem \ref{OH}.
\end{proof}

\section{The Auslander-Buchsbaum formula}

In this section we give a necessary and sufficient condition for
our homological flat dimensions to satisfy a formula of
Auslander-Buchsbaum type. Our main result is Theorem \ref{AB}
below.  Recall that
$\grade(\fp,M)=\inf\{i|\Ext^i_R(R/\fp,M)\neq0\}$.

\begin{lem}\label{AB2} Let $R$ be a local ring,
and $M$  an $R$-module of finite depth. Then
$\depth_RM\le\depth_{R_{\fp}}M_{\fp}+\dim R/\fp$ for all
$\fp\in\Spec(R)$ if and only if $\depth_RM\le\grade(\fp,M)+\dim
R/\fp$ for all $\fp\in\Spec(R)$.
\end{lem}
\begin{proof} The only if part is trivial. For the if part we show
that for modules $M$ and $N$ with $N$ finite we have
$\Ext^i_R(N,M)=0$ for $i<\depth_RM-\dim N$. We do this by
induction on $\dim N$. If $\dim N=0$, then $N$ has finite length,
and in this case an easy induction proves the result. Now let
$\dim N=t$. By a method similar to that of \cite[(17.1)]{M} it is
sufficient to take $N=R/\fp$ such that $\dim R/\fp=t$. Let
$s<\depth_RM-t\le\depth_{R_{\fp}}M_{\fp}$. We have to show that
$E=\Ext^s_R(R/\fp,M)=0$. If $E\neq0$, there is a non-zero element
$e\in E$. Since $E_{\fp}=0$ there is an element $u\in
R\backslash\fp$ such that $ue=0$. Now the exact sequence
$0\rightarrow R/\fp\stackrel{.u}{\to}R/\fp\rightarrow
N'\rightarrow 0$, gives rise the exact sequence
$$\Ext^s_R(N',M)\rightarrow\Ext^s_R(R/\fp,M)\stackrel{.u}{\to}\Ext^s_R(R/\fp,M),$$
in which the left most module equal to zero by the induction
hypothesis. So $u$ is injective and therefore $e=0$, which is a
contradiction.
\end{proof}
\begin{prop}\label{AB3} Let $M$ be an $R$-module such that $\rfd_RM+\depth_RM=\depth R$. Then
$$\depth_RM\le\grade(\fp,M)+\dim R/\fp$$ for all
$\fp\in\Spec(R)$. The converse is true over Cohen-Macaulay rings.
\end{prop}
\begin{proof} Let $\fp\in\Spec(R)$ be an arbitrary prime ideal.
Therefore we have
$$\depth
R_{\fp}-\depth_{R_{\fp}}M_{\fp}\le\rfd_{R_{\fp}}M_{\fp}\le\rfd_RM=\depth
R-\depth_RM.$$ So that
\begin{align*}
\depth_RM\le & \depth R-\depth R_{\fp}+\depth_{R_{\fp}}M_{\fp} \\[1ex]
         \le & \depth_{R_{\fp}}M_{\fp}+\dim R/\fp.
\end{align*}
Now from Lemma \ref{AB2} we obtain that
$\depth_RM\le\grade(\fp,M)+\dim R/\fp$ for all $\fp\in\Spec(R)$.

Next suppose that $R$ is a Cohen-Macaulay ring. Choose a prime
ideal $\fp\in\Spec(R)$ such that $\rfd_RM= \depth
R_{\fp}-\depth_{R_{\fp}}M_{\fp}$.  Then from the hypothesis and
\cite[Page 250]{M} we have:
\begin{align*}
\rfd_RM= & \depth R_{\fp}-\depth_{R_{\fp}}M_{\fp} \\[1ex]
       = & \dim R_{\fp}-\depth_{R_{\fp}}M_{\fp} \\[1ex]
       \le & \dim R-\dim R/\fp-\grade(\fp,M) \\[1ex]
       \le & \dim R-\depth_RM \\[1ex]
       = & \depth R-\depth_RM \\[1ex]
       \le & \rfd_RM,
\end{align*}
which completes the proof.
\end{proof}

Combining Proposition \ref{AB3}, and \cite[Theorem (3.4)]{CFF} we
have

\begin{cor} Let $R$ be a ring and $M$ be an $R$-module. The
following then are equivalent:
\begin{itemize}
\item[(i)] $\rfd_RM+\depth_RM=\depth R$ for every $R$-module $M$ of finite depth.

\item[(ii)] $R$ is a Cohen-Macaulay ring and $\depth_RM\le\grade(\fp,M)+\dim R/\fp$ for
every $R$-module $M$ of finite depth, and for all $\fp\in\Spec(R)$.
\end{itemize}
\end{cor}

\noindent This is an extension of \cite[Theorem (3.4)]{CFF}.

\noindent Now we state the main result of this section.

\begin{thm}\label{AB} Let $R$ be a Cohen-Macaulay local ring and let $M$ be an
$R$-module of finite $\H\fd_RM$ for $\H=\mbox{CI}$, $\mbox{G}^*$,
$\CM^*$, and $\CM$. Then $\H\fd_RM+\depth_RM=\depth R$ if and only
if $\depth_RM\le\grade(\fp,M)+\dim R/\fp$ for all $\fp\in\Supp(M)$.
\end{thm}
Dual to the Proposition \ref{AB3} one can prove the following.

\begin{prop} Let $R$ be a local ring,
and $M$  an $R$-module such that $\mbox{Ch}\id_RM+\depth_RM=\depth
R$. Then
$$\width_RM\le\width_{R_{\fp}}M_{\fp}+\dim R/\fp$$ for all
$\fp\in\Spec(R)$. The converse is true over Cohen-Macaulay rings.
\end{prop}

\begin{center} {\bf ACKNOWLEDGMENT}
\end{center}
\noindent The authors would like to thank H. B. Foxby and S.
Sather-Wagstaff for their useful comments, and Mehrdad Shahshahani
for his careful reading of this paper.

\end{document}